\documentclass[11pt]{amsart}
\usepackage{epsfig,amsmath,amssymb, xspace, graphicx, url, amscd, euscript, mathrsfs,color}
\input xy
\xyoption{all}
\usepackage{tikz}
\usetikzlibrary{intersections}
\usetikzlibrary{arrows}
\usetikzlibrary{decorations.pathmorphing}

\newtheorem{theorem}[subsection]{Theorem}
\newtheorem{proposition}[subsection]{Proposition}
\newtheorem{proposition-definition}[subsection]
{Proposition-Definition}
\newtheorem{corollary}[subsection]{Corollary}

\theoremstyle{definition}
\newtheorem{definition}[subsection]{Definition}

\newtheorem{example}[subsection]{Example}

\newtheorem{remark}[subsection]{Remark}

\newtheorem*{ack}{Acknowledgements}

\newcommand{\scr}[1]{\mathbf{\EuScript{#1}}}

\newcommand\cB{\mathcal{B}}
\newcommand\cC{\mathcal{C}}

\newcommand\cK{\mathcal{K}}
\newcommand\cL{\mathcal{L}}
\newcommand\cM{\mathcal{M}}

\newcommand\cO{\mathcal{O}}
\newcommand\cP{\mathcal{P}}

\newcommand\cX{\mathcal{X}}

\renewcommand\frm{\mathfrak{m}}

\begin{document}

\title{A construction of generalized elliptic curves from twisted curves}
\author[A.~Niles]{Andrew Niles}
\address{Department of Mathematics \& Computer Science\\
College of the Holy Cross\\
1 College Street\\
Worcester, MA 01610\\
USA}
\email{aniles@holycross.edu}
\date{\today}

\begin{abstract}
We introduce a procedure for constructing a generalized elliptic curve from a genus-one twisted curve, and we use this procedure to define an explicit, modular isomorphism between two compactifications of moduli stacks of elliptic curves with level structure.
\end{abstract}

\maketitle

\section{Introduction}

Modular curves and their compactifications play a central role in modern number theory. The theory of moduli spaces of elliptic curves with level structure, in characteristics not dividing the level, was put on a sound footing by Igusa in the papers \cite{I1}, \cite{I2}, \cite{I3}, \cite{I4}. The paper \cite{DR} develops this theory in terms of algebraic stacks, and gives a moduli interpretation of the cusps of the compactified stacks (in terms of level structures on generalized elliptic curves), again in characteristics not dividing the level.

Level structures on elliptic curves were reformulated as Drinfeld level structures in \cite{KM1}, where the authors constructed normal, regular, separated algebraic stacks over $\mathrm{Spec}(\mathbb{Z})$ classifying elliptic curves equipped with Drinfeld level structure. The cusps of these stacks were interpreted in \cite{C} as classifying Drinfeld level structures on generalized elliptic curves. We review this theory briefly in \S\ref{prelim} below.

The moduli stacks of elliptic curves with Drinfeld structure in \cite{KM1} arise as locally closed substacks of certain moduli stacks of twisted stable maps; this is mentioned briefly in \cite[6.9]{AOV2} and is examined in detail in \cite{N}. These stacks of twisted stable maps are complete, so this gives an alternate way of compactifying moduli stacks of elliptic curves with Drinfeld structure. For example, the moduli stack $\scr Y_1(N)$, which classifies elliptic curves equipped with $[\Gamma_1(N)]$-structure, has two compactifications: the stack $\scr X_1(N)$ considered in \cite{C}, and its closure $\scr X_1^{\mathrm{tw}}(N)$ in a certain moduli stack of twisted stable maps (see \cite[4.6]{N}).

In \cite{N}, we constructed an explicit isomorphism $\scr X_1(N) \rightarrow \scr X_1^{\mathrm{tw}}(N)$, with a natural moduli interpretation in terms of the objects classified by these two stacks. The stack $\scr X_1^{\mathrm{tw}}(N)$ classifies pairs $(\cC/S,\phi)$, where $\cC/S$ is a genus-$1$ twisted curve and $\phi$ is a Drinfeld $[\Gamma_1(N)]$-structure on $\mathrm{Pic}^0_{\cC/S}$ (see \cite[4.5]{N}). In \cite{N} we exhibited how, given a generalized elliptic curve $E/S$ and a Drinfeld $[\Gamma_1(N)]$-structure $P$ on $E$, one may naturally produce a genus-$1$ twisted curve $\cC/S$ and a Drinfeld $[\Gamma_1(N)]$-structure $\phi$ on $\mathrm{Pic}^0_{\cC/S}$.

The main result of this paper, occupying all of \S\ref{construction}, is a procedure for canonically producing generalized elliptic curves from genus-$1$ twisted curves. This construction is inverse to the (more straightforward) procedure for producing genus-$1$ twisted curves from generalized elliptic curves which we employed in \cite{N}. In particular, in \S\ref{apply} we use this construction to give a moduli interpretation for the inverse to the isomorphism $\scr X_1(N) \rightarrow \scr X_1^{\mathrm{tw}}(N)$ constructed in \cite{N}, though we hope the construction of \S\ref{construction} will prove to be of independent interest.

In \S\ref{general} we examine how our construction might extend to twisted curves of higher genus, and we discuss some complications which arise. We hope to study the higher-genus case in future work.

\begin{ack}
This paper builds on the results of the author's doctoral dissertation, and the author received many helpful comments and suggestions from his doctoral advisor, Martin Olsson, during the early stages of this project.
\end{ack}

\section{Preliminaries}\label{prelim}

\subsection*{Generalized elliptic curves and Drinfeld level structures}
We briefly recall the notion of a generalized elliptic curve and the key properties of moduli stacks of Drinfeld level structures on generalized elliptic curves; for more details the reader should consult \cite{DR}, \cite{KM1}, \cite{C} or \cite[\S2]{N}.

\begin{definition}
A \textit{Deligne-Rapoport (DR) semistable curve of genus $1$} over a scheme $S$ is a proper, flat, finitely presented morphism of schemes $f:C \rightarrow S$, all of whose geometric fibers are non-empty, connected semistable curves with trivial dualizing sheaves.
\end{definition}

The geometric fibers of a DR semistable curve of genus $1$ are either smooth genus-$1$ curves or N\'eron polygons. Over a base scheme $S$, the \textit{standard N\'eron $N$-gon} $C_N/S$ (for any $N\geq 1$) is obtained from $\widetilde{C}_N := \mathbb{P}^1_S \times \mathbb{Z}/(N)$ by ``gluing'' the $0$-section in the $i^{\mathrm{th}}$ copy of $\mathbb{P}^1_S$ to the $\infty$-section in the $(i+1)^{\mathrm{th}}$ copy of $\mathbb{P}^1_S$:

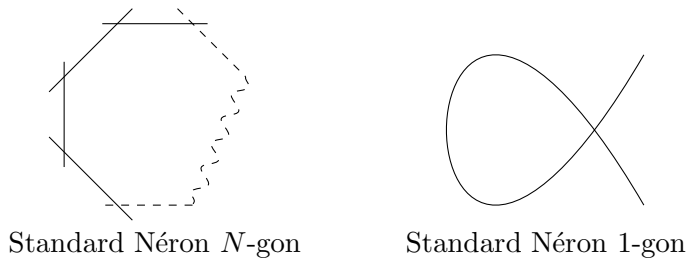
\begin{center}
\begin{tikzpicture}
\clip (-1.5,-1) rectangle (7.6,3);
\draw (0.2,-0.2) -- (-0.907,0.907);
\draw (-0.707,0.507) -- (-0.707,1.907);
\draw (-0.907,1.507) -- (0.2,2.614);
\draw (-0.2,2.414) -- (1.2,2.414);
\begin{scope}[dashed]
\draw (0.8,2.614) -- (1.707,1.707);
\draw (1,0) -- (-0.2,0);
\draw [decorate,decoration=snake] (1.707,1.707) -- (1,0);
\end{scope}
\draw [name path=curve] (7,2) .. controls (3.5,-4) and (3.5,6) .. (7,0);
\node at (0.5,-0.2) [inner sep=0pt,label=270:Standard N\'eron $N$-gon] {};
\node at (5.7,-0.2) [inner sep=0pt,label=270:Standard N\'eron $1$-gon] {};
\end{tikzpicture}
\end{center}

The natural multiplication action of $\mathbb{G}_m$ on $\mathbb{P}^1_S$, together with the action of $\mathbb{Z}/(N)$ on itself via its group law, determines an action of the group scheme $\mathbb{G}_m \times \mathbb{Z}/(N)$ on $\mathbb{P}^1_S \times \mathbb{Z}/(N)$, which descends uniquely to an action of $\mathbb{G}_m \times \mathbb{Z}/(N) = C_N^{\mathrm{sm}}$ on $C_N$ (\cite[II.1.9]{DR}).

\begin{definition}
A \textit{generalized elliptic curve} over a scheme $S$ is a genus-$1$ DR semistable curve $E/S$, equipped with a morphism $E^{\mathrm{sm}} \times E \rightarrow E$ and a section $0_E \in E^{\mathrm{sm}}(S)$ such that the restriction $E^{\mathrm{sm}} \times E^{\mathrm{sm}} \rightarrow E^{\mathrm{sm}}$ makes $E^{\mathrm{sm}}$ a commutative group scheme over $S$ with identity $0_E$, and such that on any singular geometric fiber $E_{\overline{s}}$, translation by any rational point of $E^{\mathrm{sm}}_{\overline{s}}$ acts via a rotation on the graph $\Gamma(E_{\overline{s}})$ (\cite[I.3.5]{DR}) of the irreducible components of $E_{\overline{s}}$.
\end{definition}

By \cite[II.1.15]{DR}, a generalized elliptic curve over an algebraically closed field is either a smooth elliptic curve or a N\'eron $N$-gon (for some $N\geq 1$) with the action described above.

Recall that if $X$ is a scheme over $S$, a \textit{relative effective Cartier divisor} in $X$ over $S$ is an effective Cartier divisor in $X$ which is flat over $S$. The following is the reformulation (as a Drinfeld level structure) of the na\"ive notion of a ``point of exact order $N$'' on a generalized elliptic curve, introduced in \cite{KM1} for smooth elliptic curves and in \cite{C} for generalized elliptic curves.

\begin{definition}\label{gamma1elliptic}
Let $E/S$ be a generalized elliptic curve. A \textit{$[\Gamma_1(N)]$-structure} on $E$ is a section $P \in E^{\mathrm{sm}}(S)$ such that:
\begin{itemize}
  \item $N\cdot P = 0_E$; 
  \item the relative effective Cartier divisor 
  \begin{displaymath}
  D := \sum_{a \in \mathbb{Z}/(N)} [a\cdot P]
  \end{displaymath}
  in $E^{\mathrm{sm}}$ is a subgroup scheme; and 
  \item for every geometric point $\overline{p} \rightarrow S$, $D_{\overline{p}}$ meets every irreducible component of $E_{\overline{p}}$.
\end{itemize}
We write $\scr X_1(N)$ for the stack over $\mathrm{Spec}(\mathbb{Z})$ associating to a scheme $S$ the groupoid of pairs $(E,P)$, where $E/S$ is a generalized elliptic curve and $P$ is a $[\Gamma_1(N)]$-structure on $E$. We write $\scr Y_1(N)$ for the substack classifying such pairs where $E/S$ is a smooth elliptic curve.
\end{definition}

\begin{theorem}[$\textrm{\cite[3.1.7, 3.2.7, 3.3.1, 4.1.1]{C}}$]
$\scr X_1(N)$ is a regular Deligne-Mumford stack, proper and flat over $\mathrm{Spec}(\mathbb{Z})$ of pure relative dimension $1$.
\end{theorem}

In particular, it follows (see \cite[4.1.5]{C}) that $\scr X_1(N)$ is canonically identified with the normalization of $\overline{\cM}_{1,1}$ in the normal Deligne-Mumford stack $\scr X_1(N)|_{\mathbb{Z}[1/N]}$, as in \cite[\S8.6]{KM1} and \cite[\S$\textrm{VII.2}$]{DR}.

Similarly, one may reformulate the notion of other na\"ive level structures on (generalized) elliptic curves (e.g. full level structures, or cyclic order-$N$ subgroups) as Drinfeld level structures, leading to other regular algebraic stacks (such as $\scr X(N)$ and $\scr X_0(N)$) which are also proper and flat over $\mathrm{Spec}(\mathbb{Z})$ of pure relative dimension $1$. We will focus on $\scr X_1(N)$ in this paper, but our arguments naturally and immediately extend to these other moduli stacks as well; see \cite[$\textrm{\S5, \S7}$]{N}.

\subsection*{Twisted curves}
Twisted curves are introduced in \cite{AOV2} in order to construct proper moduli spaces of stable maps to a fixed stack. If one were to only consider morphisms from nodal curves, the resulting moduli stacks would be separated but not proper; allowing the source curves to acquire some stacky structure solves the problem, producing complete moduli stacks of \textit{twisted stable maps}. We will not review the full theory of twisted stable maps here; we refer the reader to \cite{AV} and \cite{AOV2} for the original development of the theory, or to \cite[\S3]{N} for a summary.

The notion of tameness used in the following definition is that introduced in \cite{AOV1}.

\begin{definition}[$\textrm{\cite[\S2]{AOV2}}$]\label{twistedcurvedef}
An \textit{$n$-marked twisted curve} over a scheme $S$ is a proper tame stack $\cC$ over $S$, with connected geometric fibers all of dimension $1$, and with coarse space $f:C \rightarrow S$ a nodal curve over $S$; together with $n$ closed substacks 
\begin{displaymath}
\{\Sigma_i \subset \cC\}_{i = 1}^n 
\end{displaymath}
which are fppf gerbes over $S$ mapping to $n$ markings $\{p_i \in C^{\mathrm{sm}}(S)\}$, such that:
\begin{itemize}
  \item the preimage in $\cC$ of the complement $C' \subset C$ of the markings and singular locus of $C/S$ maps isomorphically onto $C'$; 
  \item if $\overline{p} \rightarrow C$ is a geometric point mapping to the image in $C$ of a marking $\Sigma_i \subset \cC$, then 
  \begin{displaymath}
  \mathrm{Spec}(\cO_{C,\overline{p}}) \times_C \cC \simeq [D^{\mathrm{sh}}/\mu_r]
  \end{displaymath}
  for some $r\geq1$, where $D^{\mathrm{sh}}$ is the strict Henselization at $(\frm_{S,f(\overline{p})},z)$ of 
  \begin{displaymath}
  D = \mathrm{Spec}(\cO_{S,f(\overline{p})}[z]) 
  \end{displaymath}
  and $\zeta\in \mu_r$ acts by $z\mapsto \zeta\cdot z$; and 
  \item if $\overline{p} \rightarrow C$ is a geometric point mapping to a node of $C$, then
  \begin{displaymath}
  \mathrm{Spec}(\cO_{C,\overline{p}})\times_C \cC \simeq [D^{\mathrm{sh}}/\mu_r]
  \end{displaymath}
  for some $r\geq 1$, where $D^{\mathrm{sh}}$ is the strict Henselization at $(\frm_{S,f(\overline{p})},x,y)$ of 
  \begin{displaymath}
  D = \mathrm{Spec}(\cO_{S,f(\overline{p})}[x,y]/(xy-t))
  \end{displaymath}
  for some $t \in \frm_{S,f(\overline{p})}$, and $\zeta\in \mu_r$ acts by $x \mapsto \zeta\cdot x$ and $y\mapsto \zeta^{-1}\cdot y$.
\end{itemize}
We say a twisted curve $\cC/S$ has \textit{genus $g$} if the geometric fibers of its coarse space $C/S$ have arithmetic genus $g$, and we say an $n$-marked genus-$g$ twisted curve $\cC/S$ is \textit{stable} if the genus-$g$ curve $C/S$ with the markings $\{p_i\}$ is an $n$-marked genus-$g$ stable curve over $S$.
\end{definition}

\begin{example}\label{stacky1gon}
Over any base scheme $S$, consider a N\'eron $1$-gon $C/S$ as above. Then $C^{\mathrm{sm}} \cong \mathbb{G}_m$, and $C$ admits the structure of a generalized elliptic curve, with an action $C^{\mathrm{sm}} \times C \rightarrow C$ extending the group scheme structure of $\mathbb{G}_m$. For any positive integer $N$, the inclusion $\mu_N \subset \mathbb{G}_m$ determines an action of $\mu_N$ on $C$, and the stack quotient $\cC := [C/\mu_N]$ is a genus-$1$ twisted stable curve over $S$, with coarse space $f:C'\rightarrow S$ a N\'eron $1$-gon. If $\overline{p} \rightarrow C'$ is a geometric point mapping to a node of $C'$, then 
\begin{displaymath}
\mathrm{Spec}(\cO_{C',\overline{p}}) \times_{C'} \cC \simeq [D^{\mathrm{sh}}/\mu_N],
\end{displaymath}
where $D^{\mathrm{sh}}$ denotes the strict Henselization of 
\begin{displaymath}
D := \mathrm{Spec}(\cO_{S,f(\overline{p})}[x,y]/(xy))
\end{displaymath}
at the point $(\mathfrak{m}_{S,f(\overline{p})},x,y)$ and $\zeta \in \mu_N$ acts by $x \mapsto \zeta \cdot x$ and $y \mapsto \zeta^{-1} \cdot y$. As in \cite{N}, we will refer to this twisted curve as the \textit{standard $\mu_N$-stacky N\'eron $1$-gon over $S$}.
\end{example}

\begin{center}
\begin{tikzpicture}
\clip (-1,-1.6) rectangle (5,1.2);
\draw [name path=curve] (3,1) .. controls (-0.5,-5) and (-0.5,5) .. (3,-1);
\node at (1.8,-1) [inner sep=0pt,label=270: Standard $\mu_N$-stacky N\'eron $1$-gon] {};
\node at (2.35,0) [circle,draw,fill,inner sep=2pt,label=0:$\mu_{N}$] {};
\end{tikzpicture}
\end{center}

Given a twisted curve $\cC$ over a scheme $S$, we write $\scr Pic_{\cC/S}$ for the stack associating to $T/S$ the groupoid of line bundles on $\cC\times_S T$.

\begin{proposition}[$\textrm{\cite[2.7]{AOV2}}$]
$\scr Pic_{\cC/S}$ is a smooth, locally finitely presented algebraic stack over $S$. For any $\cL \in \scr Pic_{\cC/S}(T)$, the group scheme $\mathrm{Aut}_T(\cL)$ is canonically isomorphic to $\mathbb{G}_{m,T}$. 
\end{proposition}

Write $\mathrm{Pic}_{\cC/S}$ for the rigidification of this stack with respect to $\mathbb{G}_m$, in the sense of \cite[\S5.1]{ACV} and \cite[Appendix A]{AOV1}; $\mathrm{Pic}_{\cC/S}$ is the relative Picard functor of $\cC/S$. From the analysis of $\mathrm{Pic}_{\cC/S}$ in \cite[\S2]{AOV2} we have:

\begin{proposition}\label{picardexact}
$\mathrm{Pic}_{\cC/S}$ is a smooth group scheme over $S$, and if $\pi:\cC \rightarrow C$ is the coarse space of $\cC/S$, then there is a short exact sequence of group schemes over $S$
\begin{displaymath}
0 \rightarrow \mathrm{Pic}_{C/S} \stackrel{\pi^*}{\rightarrow} \mathrm{Pic}_{\cC/S} \rightarrow W \rightarrow 0,
\end{displaymath}
with $W$ quasi-finite and \'etale over $S$. 
\end{proposition}

In fact, as explained in \cite[3.14]{N}, $W$ is the sheaf associated to the presheaf $T \mapsto H^0(C_T,\mathbf{R}^1\pi_* \mathbb{G}_m)$ (where we still write $\pi: \cC_T \rightarrow C_T$ for the morphism induced by base change from $\pi:\cC \rightarrow C$).

For any integer $N$ annihilating $W$, we have a natural morphism 
\begin{displaymath}
\times N: \mathrm{Pic}_{\cC/S} \rightarrow \mathrm{Pic}_{C/S}.
\end{displaymath}

\begin{definition}[$\textrm{\cite[2.11]{AOV2}}$]
The \textit{generalized Jacobian} of $\cC$ is 
\begin{displaymath}
\mathrm{Pic}^0_{\cC/S} := \mathrm{Pic}_{\cC/S} \times_{ \times N, \mathrm{Pic}_{C/S}} \mathrm{Pic}^0_{C/S}, 
\end{displaymath}
where $\mathrm{Pic}^0_{C/S}$ is the fiberwise connected component of the identity in the group scheme $\mathrm{Pic}_{C/S}$.
\end{definition}

$\mathrm{Pic}^0_{\cC/S}$ is independent of $N$, and we view it as classifying line bundles of fiberwise degree $0$ on $\cC/S$.

\begin{example}
If $\cC/S$ is a standard $\mu_N$-stacky N\'eron $1$-gon, then $\mathrm{Pic}^0_{\cC/S} \simeq \mathbb{G}_m \times \mathbb{Z}/(N)$; see \cite[3.16]{N}.
\end{example}

\section{Earlier results}

The main result of \cite{AOV2} is the construction, for any finitely presented, tame, proper algebraic stack $\cX$ with finite inertia, a complete moduli stack $\cK_{g,n}(\cX)$ classifying twisted stable maps from $n$-marked, genus-$g$ twisted curves to $\cX$. We do not require these stacks in full generality here, but we will consider one closely related stack. The stack $\overline{\cK}_{1,1}(\cB \mu_N)$ is the rigidification (in the sense of \cite[\S5.1]{ACV} and \cite[Appendix A]{AOV1}) of $\cK_{1,1}(\cB \mu_N)$ with respect to the copy of $\mu_N$ which canonically lies in the center of the group scheme of automorphisms of every object of $\cK_{1,1}(\cB \mu_N)$. It has a concrete interpretation as follows (see \cite[3.20]{N}), which for the purposes of this paper we use as our definition:

\begin{definition}
$\overline{\cK}_{1,1}(\cB \mu_N)$ is the stack over $\mathrm{Spec}(\mathbb{Z})$, associating to a scheme $S$ the groupoid of pairs $(\cC,\phi)$, where $\cC/S$ is a $1$-marked, genus-$1$ twisted stable curve, and $\phi: \mathbb{Z}/(N) \rightarrow \mathrm{Pic}^0_{\cC/S}$ is a group scheme homomorphism, such that the relative effective Cartier divisor 
\begin{displaymath}
\sum_{a \in \mathbb{Z}/(N)} [\phi(a)]
\end{displaymath}
meets every component of each geometric fiber of $\mathrm{Pic}^0_{\cC/S} \rightarrow S$.
\end{definition}

If $E/S$ is an elliptic curve, then there is a canonical isomorphism $E \simeq \mathrm{Pic}^0_{E/S}$. If $P$ is a $[\Gamma_1(N)]$-structure on $E$, we may therefore view it as a $[\Gamma_1(N)]$-structure on $\mathrm{Pic}^0_{E/S}$. This defines a group scheme homomorphism 
\begin{align*}
\mathbb{Z}/(N) & \rightarrow \mathrm{Pic}^0_{E/S} \\
1 & \mapsto P ,
\end{align*}
which in turn defines a canonical morphism of algebraic stacks $\scr Y_1(N) \hookrightarrow \overline{\cK}_{1,1} (\cB \mu_N)$. This is easily seen to be an immersion.

\begin{definition}[$\textrm{\cite[4.5]{N}}$]
Let $\cC/S$ be a $1$-marked genus-$1$ twisted stable curve with no stacky structure at its marking. A \textit{$[\Gamma_1(N)]$-structure} on $\cC$ is a group scheme homomorphism $\phi: \mathbb{Z}/(N) \rightarrow \mathrm{Pic}^0_{\cC/S}$ such that:
\begin{itemize}
  \item the relative effective Cartier divisor 
  \begin{displaymath}
  D := \sum_{a \in \mathbb{Z}/(N)} [\phi(a)]
  \end{displaymath}
  in $\mathrm{Pic}^0_{\cC/S}$ is an $S$-subgroup scheme, and 
  \item for every geometric point $\overline{p} \rightarrow S$, $D_{\overline{p}}$ meets every irreducible component of $(\mathrm{Pic}^0_{\cC/S})_{\overline{p}} = \mathrm{Pic}^0_{\cC_{\overline{p}}/k(\overline{p})}$.
\end{itemize}
We define $\scr X^{\mathrm{tw}}_1(N)$ to be the stack over $\mathrm{Spec}(\mathbb{Z})$ associating to a scheme $S$ the groupoid of pairs $(\cC,\phi) \in \overline{\cK}_{1,1}(\cB\mu_N)(S)$ such that $\phi$ is a $[\Gamma_1(N)]$-structure on $\cC$.
\end{definition}

\begin{theorem}[$\textrm{\cite[4.6, \S6]{N}}$]\label{oldtheorem}
$\scr X^{\mathrm{tw}}_1(N)$ is the closure in $\overline{\cK}_{1,1} (\cB \mu_N)$ of $\scr Y_1(N)$, and there is a canonical isomorphism $\scr X_1(N) \rightarrow \scr X^{\mathrm{tw}}_1(N)$ with a natural moduli interpretation.
\end{theorem}

We briefly recall the strategy in constructing the isomorphism $\scr X_1(N) \rightarrow \scr X^{\mathrm{tw}}_1(N)$. If $E/S$ is a generalized elliptic curve, all of whose geometric fibers are either smooth elliptic curves or N\'eron $d$-gons (for a fixed $d$ dividing $N$), then the stack quotient $\cC := [E/E^{\mathrm{sm}}[d]]$ is a $1$-marked, genus-$1$ twisted curve over $S$; the smooth geometric fibers of $\cC/S$ are the same as those of $E$, while its singular geometric fibers are all standard $\mu_d$-stacky N\'eron $1$-gons. In \cite[\S6]{N}, the author constructed a canonical isomorphism $E^{\mathrm{sm}}[N] \cong \mathrm{Pic}^0_{\cC/S}[N]$, and the isomorphism $\scr X_1(N) \rightarrow \scr X^{\mathrm{tw}}_1(N)$ sends the generalized elliptic curve $E/S$ with a $[\Gamma_1(N)]$-structure $P$ to the twisted curve $\cC$ with the $[\Gamma_1(N)]$-structure induced by $P$ via the isomorphism $E^{\mathrm{sm}}[N] \cong \mathrm{Pic}^0_{\cC/S}[N]$.

\section{Constructing generalized elliptic curves from twisted curves}\label{construction}

This section is devoted to our main result, a procedure for constructing a generalized elliptic curve from a genus-$1$ twisted curve.

Specifically, let $B = \mathrm{Spec}(A)$ be the spectrum of a complete, strictly Henselian local ring $A$, and let $\cC/B$ be a $1$-marked genus-$1$ twisted curve, with no stacky structure at its marking. To $\cC$ we will associate a canonical generalized elliptic curve $E/B$. This is only nontrivial if the closed fiber of the coarse space $f: C \rightarrow B$ is singular (a smooth genus-$1$ twisted curve with non-stacky marking is simply an elliptic curve), so we will assume for the remainder of the proof that the closed fiber of the coarse space $C/B$ is a genus-$1$ rational nodal curve. We will employ a slightly modified version of a construction appearing in \cite{O} and \cite{A2}.

For each $n$, let $A_n = A/\mathfrak{m}^{n+1}$, $B_n = \mathrm{Spec}(A_n)$ and $\cC_n = \cC \times_B B_n$. Since $\cC$ is a twisted curve, there is a unique integer $d$ such that, if $\overline{p} \rightarrow C$ is a geometric point mapping to the node of $C_0$, we may write 
\begin{displaymath}
\mathrm{Spec}(\cO_{C,\overline{p}}) \times_C \cC \simeq [D^{\mathrm{sh}}/\mu_d]. 
\end{displaymath}
Here $D^{\mathrm{sh}}$ is the strict Henselization at $(\mathfrak{m}_{B,f(\overline{p})},x,y)$ of 
\begin{displaymath}
\mathrm{Spec}(\cO_{B,f(\overline{p})}[x,y]/(xy-t)) 
\end{displaymath}
for some appropriate $t \in \mathfrak{m}_{B,f(\overline{p})}$, and $\zeta \in \mu_d$ acts by $x \mapsto \zeta x$ and $y \mapsto \zeta^{-1} y$. More tersely, $d$ is the index of the stack $\cC$ in the sense of \cite[2.3.3]{AH}. It follows immediately from \cite[3.17]{N} that after a finite base change, $\cC_0$ is a standard $\mu_d$-stacky N\'eron $1$-gon over the residue field $A_0 = A/\mathfrak{m}$. In particular, by \cite[3.14]{N} we have a canonical short exact sequence 
\begin{displaymath}
0 \rightarrow \mu_d \rightarrow \mathrm{Pic}^0_{\cC_n/B_n}[d] \rightarrow \mathbb{Z}/(d) \rightarrow 0.
\end{displaymath}

Fix a section $Q \in \mathrm{Pic}^0_{\cC_n/B_n}(B_n)$, lying over $1 \in \mathbb{Z}/(d)$ in the above short exact sequence, such that 
\begin{align*}
\beta: \mathbb{Z}/(d) & \rightarrow \mathrm{Pic}^0_{\cC_n/B_n} \\
1 & \mapsto Q 
\end{align*}
is a $[\Gamma_1(d)]$-structure on $\cC_n$ in the sense of \cite[4.5]{N}. The above short exact sequence implies that such a $Q$ exists, at least after a finite base change. (The canonical nature of our constructions guarantees that our ultimate construction of a generalized elliptic curve $E_n/B_n$ will descend to our original $B_n$, so such a finite base change is harmless.)

After a further finite base change, we may assume that the morphism $\beta$ factors through a morphism of Picard stacks $\alpha: \mathbb{Z}/(d) \rightarrow \scr{P}ic^0_{\cC_n/B_n}$, where 
\begin{displaymath}
\scr{P}ic^0_{\cC_n/B_n} = \scr{P}ic_{\cC_n/B_n} \times_{\mathrm{Pic}_{\cC_n/B_n}} \mathrm{Pic}^0_{\cC_n/B_n} 
\end{displaymath}
for $\scr{P}ic_{\cC_n/B_n}$ the Picard stack of line bundles on $\cC_n$. $\mathrm{Pic}^0_{\cC_n/B_n}$ is defined as in \cite[2.11]{AOV2}. Define 
\begin{displaymath}
C'_n := \underline{\mathrm{Spec}}_{\cC_n} \big( \bigoplus_{m \in \mathbb{Z}/(d)} \cL_m \big),
\end{displaymath}
where $\cL_m = \alpha(m)$. Since $\alpha$ is a morphism of Picard stacks, there are canonical isomorphisms $\cL_a \otimes \cL_b \cong \cL_{a+b}$, defining the algebra structure on $\oplus \cL_m$. The $\mathbb{Z}/(d)$-grading defines a natural action of $\mu_d$ on $C'_n$, making $C'_n$ a $\mu_d$-torsor over $\cC_n$. As a $\mu_d$-torsor over $\cC_n$, $C'_n$ is independent (up to canonical isomorphism) of the choice of $\alpha$ lifting $\beta$.

Since $Q$ is a $[\Gamma_1(d)]$-structure, $\langle Q \rangle$ meets every component of $\mathrm{Pic}^0_{\cC_n/B_n}$. This implies that if $\xi: \mathrm{Spec}(k) \rightarrow \cC_n$ is a geometric point mapping to the node of $C_0$, then the action of $\mathrm{Aut}(\xi)$ on the fiber of $\cL_1$ is effective, so $C'_n$ is representable by \cite[2.3.10]{AH}. The composition 
\begin{displaymath}
C'_n \rightarrow \cC_n \rightarrow C_n 
\end{displaymath}
is a degree-$d$ cover of the genus-$1$ rational nodal curve $C_n$, totally ramified over the node of $C_n$ and unramified elsewhere, so it follows from Riemann-Hurwitz for nodal curves that $C'_n$ is naturally a genus-$1$ nodal curve over $B_n$, with closed fiber a N\'eron $1$-gon. 

Fixing a section $0_{C'_n}$ of $C'_n$ lifting the marked section $0_{\cC_n}$ of $\cC_n$ gives $C'_n$ the structure of a generalized elliptic curve. The generalized elliptic curve arising from a different choice of $0_{C'_n}$ is canonically isomorphic to $C'_n$ with the given structure of a generalized elliptic curve. (To see this, choose an isomorphism $(C'_n)^{\mathrm{sm}} \cong \mathbb{G}_m$. If the two choices of identity element correspond to $\zeta_0,\zeta_1 \in \mathbb{G}_m$, then multiplication by $\zeta_1 \zeta_0^{-1}$ gives the desired isomorphism.) And while the structure of $C'_n$ as a $\mu_d$-torsor over $\cC_n$ depends on the choice of $Q$ used to define it, the structure of $C'_n$ as a generalized elliptic curve is clearly independent of this choice.

Let $X$ be the character lattice of $\mathrm{Pic}^0_{C'_n/B_n}$; since $C'_n$ is a genus $1$ nodal curve over $B_n$, with closed fiber a N\'eron $1$-gon, $\mathrm{Pic}^0_{C'_n/B_n}$ is noncanonically isomorphic to $\mathbb{G}_m$, and so $X$ is noncanonically isomorphic to $\mathbb{Z}$. Fix an isomorphism $X \cong \mathbb{Z}$. Let 
\begin{displaymath}
\widetilde{\cP}_n \rightarrow C'_n 
\end{displaymath}
be the $X$-torsor arising from the construction of \cite[$\textrm{\S4.2}$]{A2} and refined in \cite[$\textrm{\S4.1}$]{O} (both are generalizations of the construction of the Tate curve in \cite[$\textrm{\S VII}$]{DR}). 

Explicitly, $\widetilde{\cP}_n$ is an infinite chain of copies of $\mathbb{P}^1_{B_n}$, indexed by $X$, with the $0$-section of each copy of $\mathbb{P}^1$ meeting the $\infty$-section of an adjacent copy of $\mathbb{P}^1$ at a node. If the node of $C'_n$ is locally of the form $\mathrm{Spec}(A_n[x,y]/(xy - t_n))$, the same is true for each node of $\widetilde{\cP}_n$: 
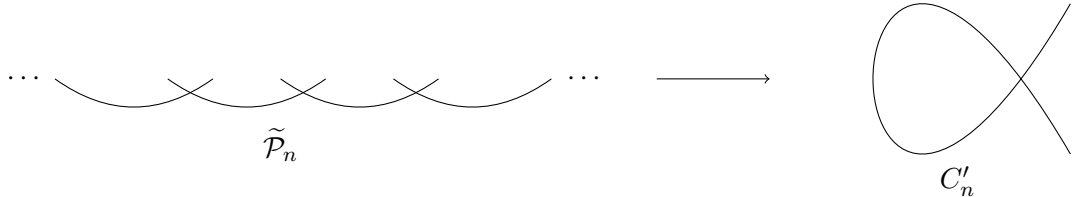
\begin{center}
\begin{tikzpicture}
\clip (-6.5,3.4) rectangle (9,6.1);
\draw [name path=curve] (-5.5,5) .. controls (-4.8,4.5) and (-4.1,4.5) .. (-3.4,5);
\draw [name path=curve] (-4,5) .. controls (-3.3,4.5) and (-2.6,4.5) .. (-1.9,5);
\draw [name path=curve] (-2.5,5) .. controls (-1.8,4.5) and (-1.1,4.5) .. (-0.4,5);
\draw [name path=curve] (-1,5) .. controls (-0.3,4.5) and (0.4,4.5) .. (1.1,5);
\draw [name path=curve] (8,6) .. controls (4.5,0) and (4.5,10) .. (8,4);
\node at (6.5,4) [inner sep=0pt,label=270:\textrm{$C'_n$}] {};
\node at (-2.5,4.5) [inner sep=0pt,label=270:$\widetilde{\cP}_n$] {};
\node at (1.15,5) [inner sep=0pt,label=0:$\cdots$] {};
\node at (-6.3,5) [inner sep=0pt,label=0:$\cdots$] {};
\draw [->] (2.5,5) -- (4,5);
\end{tikzpicture}
\end{center}
$X$ acts on $\widetilde{\cP}_n$ by translations: the action of a generator of $X \cong \mathbb{Z}$ shifts each component over to an adjacent component. Additionally the torus $\mathrm{Pic}^0_{C'_n/B_n}$ acts on $\widetilde{\cP}_n$, preserving the components. The choice of an identity element $0_{\widetilde{\cP}_n} \in \mathrm{Pic}^0_{C'_n/B_n}(B_n)$ lifting $0_{C'_n}$ determines an isomorphism $\widetilde{\cP}_n^{\mathrm{sm}} \cong \mathrm{Pic}^0_{C'_n/B_n} \times X$; so if an identity is chosen, $\widetilde{\cP}_n^{\mathrm{sm}}$ becomes a group scheme, and one may view the induced action of $\mathrm{Pic}^0_{C'_n/B_n} \times X$ on $\widetilde{\cP}_n$ as an action of $\widetilde{\cP}_n^{\mathrm{sm}}$.

Consider the quotient $E_n := \widetilde{\cP}_n/(dX)$ of $\widetilde{\cP}_n$ by the sublattice $dX \subseteq X$. $E_n$ is an $X/(dX)$-torsor over $C'_n$; under the chosen isomorphism $X \cong \mathbb{Z}$ we may view $E_n$ as a $\mathbb{Z}/(d)$-torsor over $C'_n$. Choose a section $0_{E_n} \in E_n^{\mathrm{sm}}(B_n)$ lying over $0_{C'_n}$; the choice of $0_{E_n}$ makes $E_n$ a generalized elliptic curve over $B_n$, with closed fiber a N\'eron $d$-gon. The action of $E_n^{\mathrm{sm}}$ is induced by the action of $\widetilde{\cP}_n^{\mathrm{sm}}$ on $\widetilde{\cP}_n$, by fixing any identity element $0_{\widetilde{\cP}_n}$ lying over $0_{E_n}$ (the action is clearly independent of this choice).

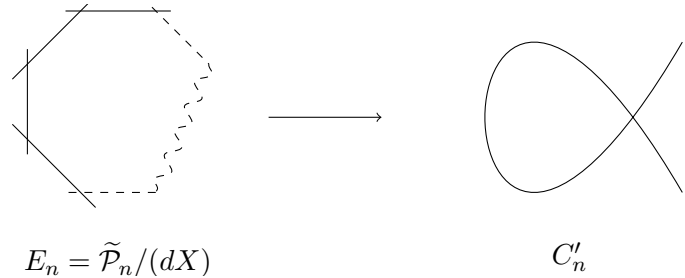
\begin{center}
\begin{tikzpicture}
\clip (-2,2.9) rectangle (9,6.5);
\draw (0.2,3.8) -- (-0.907,4.907);
\draw (-0.707,4.507) -- (-0.707,5.907);
\draw (-0.907,5.507) -- (0.2,6.614);
\draw (-0.2,6.414) -- (1.2,6.414);
\begin{scope}[dashed]
\draw (0.8,6.614) -- (1.707,5.707);
\draw (1,4) -- (-0.2,4);
\draw [decorate,decoration=snake] (1.707,5.707) -- (1,4);
\end{scope}
\draw [name path=curve] (8,6) .. controls (4.5,0) and (4.5,10) .. (8,4);
\node at (6.5,3.5) [inner sep=0pt,label=270:\textrm{$C'_n$}] {};
\node at (0.5,3.5) [inner sep=0pt,label=270:$\textrm{$E_n = \widetilde{\cP}_n/(dX)$}$] {};
\draw [->] (2.5,5) -- (4,5);
\end{tikzpicture}
\end{center}

The generalized elliptic curve arising from a different choice of $0_{E_n}$ is canonically isomorphic to $E_n$ with the given structure of a generalized elliptic curve (the canonical isomorphism being a rotation), and the same would be true had we used a different isormosphim $X \cong \mathbb{Z}$ above (the canonical isomorphism being the inversion isogeny $[-1]$, composed with a rotation if necessary). 

Compatibility of the twisted curves $\{\cC_n\}$ as $n$ varies follows from the compatibility of the generalized elliptic curves $\{C_n/B_n\}$ (which is clear from the construction) and from the compatibility of the $X$-torsors $\{\widetilde{\cP}_n \rightarrow C'_n\}$ (see \cite[\S4.1]{O}). Thus the inverse system $\{E_n/B_n\}$ algebraizes uniquely to a canonical generalized elliptic curve $E$ over $B$, completing our construction.

\begin{remark}\label{cover}
Note that by construction, $E$ comes equipped with a canonical morphism $E \rightarrow \cC$, identifying $\cC$ with the stack quotient $[E/E^{\mathrm{sm}}[d]]$.
\end{remark}

\section{An application}\label{apply}

The construction of \S\ref{construction} immediately allows us to give a moduli interpretation of the inverse to the isomorphism $\scr X_1(N) \rightarrow \scr X_1^{\mathrm{tw}}(N)$ of algebraic stacks of Theorem \ref{oldtheorem}.

To give a map $\scr X_1^{\mathrm{tw}}(N) \rightarrow \scr X_1(N)$, it suffices to define a map 
\begin{displaymath}
\scr X_1^{\mathrm{tw}}(N)(B) \rightarrow \scr X_1(N)(B), 
\end{displaymath}
in the case where $B = \mathrm{Spec}(A)$ is the spectrum of a complete, strictly Henselian local ring. 

Let 
\begin{displaymath}
(\cC,\phi) \in \scr X_1^{\mathrm{tw}}(N)(B), 
\end{displaymath}
so $\cC/B$ is a $1$-marked genus-$1$ twisted curve (with non-stacky marking) and 
\begin{displaymath}
\phi: \mathbb{Z}/(N) \rightarrow \mathrm{Pic}^0_{\cC/B} 
\end{displaymath}
is a $[\Gamma_1(N)]$-structure on $\cC$. As in \S\ref{construction}, let $B_n = \mathrm{Spec}(A/\mathfrak{m}^{n+1})$ and $\cC_n = \cC \times_B B_n$. Also let $\phi_n$ be the $[\Gamma_1(N)]$-structure on $\cC_n$ induced by $\phi$.

Let $d$ denote the index (in the sense of \cite[2.3.3]{AH}) of $\cC$, as in \S\ref{construction}. Since $d$ equals the number of components of $\mathrm{Pic}^0_{\cC_0/B_0}$, and since $\cC$ admits a $[\Gamma_1(N)]$-structure (namely $\phi$), it follows that $d$ divides $N$.

Let $E_n/B_n$ be the generalized elliptic curve associated to $\cC_n$, as constructed in \S\ref{construction}. Since $[E_n/E_n^{\mathrm{sm}}[d]] \simeq \cC_n$ canonically (see Remark \ref{cover}), the construction in the proof of \cite[6.8]{N} gives an explicit, canonical isomorphism 
\begin{displaymath}
\gamma: E_n^{\mathrm{sm}}[N] \cong \mathrm{Pic}^0_{\cC_n/B_n}[N].
\end{displaymath}
We set $P_n = \gamma^{-1}(\phi_n(1))$, where $\phi_n$ is the restriction to $\cC_n$ of the $[\Gamma_1(N)]$-structure $\phi$. $P_n$ is a $[\Gamma_1(N)]$-structure on the generalized elliptic curve $E_n$, and our construction shows that the pair $(E_n,P_n)$ is canonically determined by the pair $(\cC_n,\phi_n)$. 

Compatibility of the pairs $(\cC_n,\phi_n)$ as $n$ varies is automatic, as they are all induced by the pair $(\cC,\phi)$. This implies compatibility of the pairs $(E_n,P_n)$, and thus the inverse system $\{(E_n,P_n)\}$ algebraizes uniquely to a canonical pair 
\begin{displaymath}
(E,P) \in \scr X_1(N)(B),
\end{displaymath}
where $E$ is the generalized elliptic curve over $B$ constructed in \S\ref{construction}, and $P$ is a $[\Gamma_1(N)]$-structure on $E$. The morphism of algebraic stacks 
\begin{displaymath}
\scr X_1(N) \rightarrow \scr X_1^{\mathrm{tw}}(N)
\end{displaymath}
defined in \cite[\S6]{N} clearly sends the pair $(E_n,P_n) \in \scr X_1(N)(B_n)$ to the pair $(\cC_n,\phi_n) \in \scr X_1^{\mathrm{tw}}(N)(B_n)$ for each $n$, so it sends $(E,P) \in \scr X_1(N)(B)$ to a pair canonically isomorphic to $(\cC,\phi) \in \scr X_1^{\mathrm{tw}}(N)(B)$. Therefore we have proven:

\begin{corollary}
The morphism of algebraic stacks 
\begin{align*}
\scr X_1^{\mathrm{tw}}(N) & \rightarrow \scr X_1(N) \\
(\cC,\phi) & \mapsto (E,P) 
\end{align*}
is quasi-inverse to the isomorphism $\scr X_1(N) \rightarrow \scr X_1^{\mathrm{tw}}(N)$ defined in \cite[\S6]{N}.
\end{corollary}

We remark that the same construction of generalized elliptic curves from twisted curves gives moduli interpretations for the inverses to the isomorphisms $\scr X(N) \rightarrow \scr X^{\mathrm{tw}}(N)$ and $\scr X_0(N) \rightarrow \scr X_0^{\mathrm{tw}}(N)$ considered in \cite{N}. The argument is essentially identical to that for the inverse to $\scr X_1(N) \rightarrow \scr X_1^{\mathrm{tw}}(N)$, so we leave the details to the reader.

\section{Generalizations}\label{general}

An obvious question to ask is the extent to which the construction of \S\ref{construction} can be generalized, such as to twisted curves of higher genus. In this section we explore how one might generalize the construction, and the complications that arise in this case.

For simplicity, we work over the spectrum $S = \mathrm{Spec}(k)$ of an algebraically closed field $k$. Let $\cC/S$ be an unmarked twisted curve, whose coarse space $C/S$ is an irreducible, genus-$g$ nodal curve ($g \geq 2$) with two nodes. Suppose one node has index $d$ and the other has index $e$.

\begin{center}
\begin{tikzpicture}
\clip (-1,-1.6) rectangle (7,1.6);
\draw [name path=curve] (3,0.5) .. controls (-0.5,-5) and (-0.5,5) .. (3,-0.5);
\draw [name path=curve] (3,0.5) .. controls (3.8,1.5) and (4.6,1.5) .. (6,-0.5);
\draw [name path=curve] (3,-0.5) .. controls (3.8,-1.5) and (4.6,-1.5) .. (6,0.5);
\node at (2.67,0) [circle,draw,fill,inner sep=2pt,label=0:$\mu_{d}$] {};
\node at (5.62,0) [circle,draw,fill,inner sep=2pt,label=0:$\mu_{e}$] {};
\node at (3,-1) [inner sep=0pt,label=270: $\cC$] {};
\end{tikzpicture}
\end{center}

By \cite[3.14]{N}, we have a canonical short exact sequence 
\begin{displaymath}
0 \rightarrow \mathrm{Pic}^0_{C/S} \rightarrow \mathrm{Pic}^0_{\cC/S} \rightarrow \mathbb{Z}/(d) \times \mathbb{Z}/(e) \rightarrow 0.
\end{displaymath}

If we fix a morphism of Picard stacks $\phi: \mathbb{Z}/(d) \times \mathbb{Z}/(e) \rightarrow \scr Pic^0_{\cC/S}$ lifting a splitting of the above short exact sequence, then we get a $\mu_d \times \mu_e$-torsor 
\begin{displaymath}
C' = \underline{\mathrm{Spec}}_{\cC} \bigg( \bigoplus_{(a,b) \in \mathbb{Z}/(d) \times \mathbb{Z}/(e)} \phi(a,b) \bigg)
\end{displaymath}
over $\cC$. By the same argument as in \S\ref{construction}, $C'$ is representable.

One would then be tempted to apply the constructions of \cite[\S4.2]{A2} and \cite[\S4.1]{O} to obtain an infinite cover of a compactified Jacobian of $C'$. But this is incorrect in general, because $C'$ need not be a nodal curve at all! If it were, then the composite $C' \rightarrow \cC \rightarrow C$ would be a finite morphism of nodal curves of degree $de$; it would be unramified away from the nodes, and every point of $C'$ lying over the first (resp. second) node of $C$ would have ramification index $d$ (resp. $e$). But for arbitrary indices $d$ and $e$, an elementary Riemann-Hurwitz argument implies that no such $C'$ exists.

There is an alternative way to proceed which avoids this problem. Let $\cC/S$ be the same twisted curve considered above. Fix an isomorphism $\mathrm{Pic}^0_{C/S} \cong \mathbb{G}_m^2 \times \mathrm{Pic}^0_{C^\nu /S}$, where $C^\nu$ is the normalization of $C$. Let $X$ denote the character lattice of $\mathrm{Pic}^0_{C/S}$, so $X \cong \mathbb{Z}^2$ noncanonically. The canonical injection 
\begin{displaymath}
\mu_d \times \mu_e \hookrightarrow \mathbb{G}_m^2 \times \{1\} \subseteq \mathbb{G}_m^2 \times \mathrm{Pic}^0_{C^\nu /S} = \mathrm{Pic}^0_{C/S} 
\end{displaymath}
determines a surjection $X \twoheadrightarrow \mathbb{Z}/(d) \times \mathbb{Z}/(d)$. Let $K \subseteq X$ be the kernel of this surjection.

The constructions of \cite[\S4.2]{A2} and \cite[\S4.1]{O} provide a canonical $X$-torsor 
\begin{displaymath}
\widetilde{\cP} \rightarrow \underline{\mathrm{Jac}}_{g-1}(C/S).
\end{displaymath}
Here $\underline{\mathrm{Jac}}_{g-1}(C/S)$ is the canonical compactified Jacobian of $C$ as in \cite{A1}. We then obtain a $\mathbb{Z}/(d) \times \mathbb{Z}/(e)$-torsor 
\begin{displaymath}
\widetilde{\cP}/K \rightarrow \underline{\mathrm{Jac}}_{g-1}(C/S).
\end{displaymath}

$\widetilde{\cP}/K$ can be thought of as a higher-dimensional analogue of a generalized elliptic curve. Indeed, if we had followed this construction where $\cC$ is instead a $1$-marked genus-$1$ twisted curve with node of index $d$ as in \S\ref{construction}, $\widetilde{\cP}/K$ would be a N\'eron $d$-gon and a $\mathbb{Z}/(d)$-torsor over $C$. The primary difference between this construction and that of \S\ref{construction} is that we obtain a cover of the coarse space $C$, rather than a cover of a nodal curve $C'$ which in turn covers the twisted curve $\cC$.

On the other hand, $\widetilde{\cP}/K$ can also be thought of as a compactified Jacobian for $\cC$. Indeed, the constructions of \cite[\S4.2]{A2} and \cite[\S4.1]{O} equip $\widetilde{\cP}$ with a canonical action of $\mathrm{Pic}^0_{C/S} \times X$, and the smooth locus of $\widetilde{\cP}$ is a principal homogeneous space under this action. This induces an action of $\mathrm{Pic}^0_{C/S} \times (X/K) \cong \mathrm{Pic}^0_{\cC/S}$ on $\widetilde{\cP}/K$, and the smooth locus of $\widetilde{\cP}$ is a principal homogeneous space under this action.

This construction can clearly be extended to the case where $\cC/S$ is an unmarked twisted curve of genus $g \geq 2$, with irreducible coarse space $C/S$ and an arbitrary number of nodes. One natural question which arises is whether the requirement that $C$ be irreducible is necessary; we adopted this requirement so that the toric part of $\mathrm{Pic}^0_{C/S}$ was of maximal rank, allowing us to find a surjection $X \twoheadrightarrow \mathbb{Z}/(d) \times \mathbb{Z}/(e)$. It would also be interesting to see whether $\widetilde{\cP}/K$ has a natural moduli interpretation, analogously to that of $\underline{\mathrm{Jac}}_{g-1}(C/S)$. We do not attempt to answer these questions here, but hope to revisit them in future work.


\begin{thebibliography}{999999}

\bibitem[ACV]{ACV} D.~Abramovich, A.~Corti, A.~Vistoli, \textit{Twisted bundles and admissible covers}, Comm.~Algebra \textbf{31} (2003), no.~8, 3547-3618.

\bibitem[AH]{AH} D.~Abramovich, B.~Hassett, \textit{Stable varieties with a twist}, in \textit{Classification of Algebraic Varieties}, EMS Ser.~Congr.~Rep., Eur.~Math.~Soc.~(2011), 1-38.

\bibitem[AOV1]{AOV1} D.~Abramovich, M.~Olsson, A.~Vistoli, \textit{Tame stacks in positive characteristic}, Ann.~Inst.~Fourier (Grenoble) \textbf{58} (2008), no.~4, 1057-1091.

\bibitem[AOV2]{AOV2} D.~Abramovich, M.~Olsson, A.~Vistoli, \textit{Twisted stable maps to tame Artin stacks}, J.~Algebraic Geom.~\textbf{20} (2011), no.~2, 399-477.

\bibitem[AV]{AV} D.~Abramovich, A.~Vistoli, \textit{Compactifying the space of stable maps}, J.~Amer.~Math.~Soc.~\textbf{15} (2002), no.~1, 27-75.

\bibitem[A1]{A1} V. Alexeev, \textit{Compactified Jacobians and Torelli map}, Publ. Res. Inst. Math. Sci. \textbf{40} (2004), no. 4, 1241-1265.

\bibitem[A2]{A2} V.~Alexeev, \textit{Complete moduli in the presence of semiabelian group action}, Ann.~ of Math.~(2) \textbf{155} (2002), no.~3, 611-708.

\bibitem[Con]{C} B.~Conrad, \textit{Arithmetic moduli of generalized elliptic curves}, J.~Inst.~Math.~Jussieu \textbf{6} (2007), no.~2, 209-278.

\bibitem[DR]{DR} P.~Deligne, M.~Rapoport, \textit{Les sch\'emas de modules des courbes elliptiques}, in \textit{Modular Functions of One Variable II}, Lecture Notes in Mathematics \textbf{349}, Springer (1973), 143-316.

\bibitem[Igu1]{I1} J. Igusa, \textit{Fibre systems of Jacobian varieties. III. Fibre systems of elliptic curves.} Amer. J. Math. \textbf{81} (1959), 453-476.

\bibitem[Igu2]{I2} J. Igusa, \textit{Kroneckerian model of fields of elliptic modular functions}, Amer. J. Math. \textbf{81} (1959), 561-577.

\bibitem[Igu3]{I3} J. Igusa, \textit{On the transformation theory of elliptic functions}, Amer. J. Math. \textbf{81} (1959), 436-452.

\bibitem[Igu4]{I4} J. Igusa, \textit{On the algebraic theory of elliptic modular functions}, J. Math. Soc. Japan \textbf{20} (1968), 96-106.

\bibitem[KM]{KM1} N.~Katz, B.~Mazur, \textit{Arithmetic Moduli of Elliptic Curves}, Annals of Mathematics Studies \textbf{108}, Princeton University Press (1985).

\bibitem[Nil]{N} A.~Niles, \textit{Moduli of elliptic curves via twisted stable maps}, Algebra Number Theory \textbf{7} (2013), no.~9, 2141-2202.

\bibitem[Ols]{O} M.~Olsson, \textit{Compactifying Moduli Spaces for Abelian Varieties}, Lecture Notes in Mathematics \textbf{1958}, Springer (2008).

\end{thebibliography}
\end{document}